\documentclass[11pt,a4paper,twoside]{article}
\usepackage{amssymb}
\usepackage{amsmath}
\usepackage{fancyhdr}
\usepackage{amsthm}
\usepackage{rotate}
\usepackage{graphicx}
\usepackage{float}
\usepackage[T1]{fontenc}

\usepackage{afterpage}  %

\newtheorem{theorem}{Theorem}[section]
\newtheorem{lemma}[theorem]{Lemma}
\newtheorem{proposition}[theorem]{Proposition}

\newtheorem{example}[theorem]{Example}
\newtheorem{definition}[theorem]{Definition}
\theoremstyle{remark}

\begin{document}

\title{Rigid and super rigid quasigroups}
\author{Andriy I. Deriyenko, Ivan I. Deriyenko and Wieslaw A.
Dudek} \date{}\maketitle

{\bf Abstract.} {\footnotesize The paper deals with quasigroups
having a trivial group of automorphisms and a trivial group of
autotopisms. Examples of such quasigroups and methods of their
verification are given. }

\section{Introduction}

Let $Q=\{1,2,3,\ldots,n\}$ be a finite set, $\varphi$ and $\psi$
permutations of $Q$. The multiplication (composition) of
permutations is defined as $\varphi\psi(x)=\varphi(\psi(x))$.
Permutations will be written in the form of cycles and cycles will
be separated by points, e.g.

$$
\varphi=\left(
\begin{array}{cccccc}
1 & 2 & 3 & 4 & 5 & 6\\
3 & 1 & 2 & 5 & 4 & 6
\end{array}
\right)=(123.45.6.)
$$

As it is well known, any permutation $\varphi$ of the set $Q$ of
order $n$ can be decomposed into $r\leqslant n$ cycles of the
length $k_1,k_2,\ldots,k_r$ and $k_1 + k_2 + \ldots + k_r = n$. We
denote this fact by
$$
Z(\varphi)=[k_1,k_2,\ldots,k_r].
$$
Two permutations are conjugate if and only if they have the same
number of cycles of each length (Theorem 5.1.3 in  \cite{4}). So,
for any two permutations $\varphi$ and $\psi$ we have
$$
Z(\varphi)=Z(\psi)\longleftrightarrow \beta\varphi\beta^{-1}=\psi.
$$
From the proof of Theorem 5.1.3 and Lemma 5.1.1 in \cite{4}
follows a method of determination of $\beta$. This method is also
used here, so let us recall it.

If $\beta\varphi\beta^{-1}=\psi$ and
$$
\begin{array}{l}
\varphi=(a_{11}\thinspace a_{12}\ldots a_{1k_1})\ldots(a_{r1}\ldots a_{rk_r}) \\[3pt]
\psi=(b_{11}\thinspace b_{12}\ldots b_{1k_1})\ldots(b_{r1}\ldots b_{rk_r})
\end{array}
$$
then, according to \cite{4}, $\beta$ has the form
\begin{equation}
\label{2}
\beta=\left(
\begin{array}{cccccccc}
a_{11} & a_{12} & \ldots & a_{1k_1} & \ldots & a_{r1} & \ldots & a_{rk_r} \\
b_{11} & b_{12} & \ldots & b_{1k_1} & \ldots & b_{r1} & \ldots & b_{rk_r}
\end{array}
\right),
\end{equation}
where the first row contains all elements of $\varphi$, the second
-- elements of $\psi$ written in the same order as in
decompositions into cycles. Replacing in $\varphi$ the cycle
$(a_{11}\thinspace a_{12}\ldots a_{1k_1})$ by $(a_{12}\thinspace
a_{13}\ldots a_{1k_1}\thinspace a_{11})$ we save the permutation
$\varphi$ but we obtain a new $\beta$. Similarly for an
arbitrary cycle of $\varphi$ and $\psi$. One can prove that in
this way we obtain all $\beta$ satisfying the equality
$\beta\varphi\beta^{-1}=\psi$.

\begin{definition}\rm Let $Q(\cdot)$ be a quasigroup. Each permutation $\varphi_i$ of
$Q$ satisfying the identity
\begin{equation}\label{t}
x\cdot \varphi_i(x)=i,
\end{equation}
where $i\in Q$, is called a
{\it track\,} or a {\it right middle translation}.
\end{definition}

Such permutations were firstly studied by V. D. Belousov
\cite{Bel} in connection with some groups associated with
quasigroups. The investigations of such permutations were
continued, for example, in \cite{3, 2} and \cite{Sch}.

The above condition says that in a Latin square $n\times n$
associated with a quasigroup $Q(\cdot)$ of order $n$ we select $n$
cells, one in each row, one in each column, containing the same
fixed element $i$. $\varphi_i(x)$ means that to find in the row
$x$ the cell containing $i$ we must select the column
$\varphi_i(x)$. It is clear that for a quasigroup $Q(\cdot)$ of
order $n$ the set $\{\varphi_1,\varphi_2,\ldots,\varphi_n\}$
uniquely determines its Latin square, and conversely, any Latin
square $n\times n$ uniquely determines the set
$\{\varphi_1,\varphi_2,\ldots,\varphi_n\}$.

Connections between tracks of isotopic quasigroups are described
in  \cite{3} and \cite{2}. Namely, if
$\{\varphi_1,\varphi_2,\ldots\varphi_n\}$ are tracks of
$Q(\cdot)$, $\{\psi_1,\psi_2,\ldots\psi_n\}$ -- tracks of
$Q(\circ)$ and
$$
\gamma(x\circ y)=\alpha(x)\cdot\beta(y), $$
then
\begin{equation}
\label{3} \varphi_{\gamma(i)}=\beta\psi_i \alpha^{-1}.
\end{equation}
So, tracks of isomorphic quasigroups ($\alpha=\beta=\gamma$) are
connected by the formula
$$
\varphi_{\alpha(i)}=\alpha\psi_{i}\alpha^{-1}.
$$
Thus, for any automorphism $\alpha$ of a quasigroup $Q(\cdot)$ we
have
\begin{equation}
\label{4}
\varphi_{\alpha(i)}=\alpha\varphi_{i}\alpha^{-1}
\end{equation}
and
\begin{equation}
\label{5} Z(\varphi_{i})=Z(\varphi_{\alpha(i)}).
\end{equation}

\begin{definition}\rm A track $\varphi_{k}$ of $Q(\cdot)$ is called {\it
special} if $Z(\varphi_{k})\neq Z(\varphi_{i})$ for all $i\in Q$,
$i\neq k$.
\end{definition}
\begin{example}\rm\label{ex1} Consider two isotopic quasigroups:
{\footnotesize{\[
\begin{array}{lcr}
\begin{array}{c|ccc}
\cdot & 1 & 2 & 3 \\ \hline\rule{0mm}{3.5mm}
1       & 1 & 2 & 3 \\
2       & 2 & 3 & 1 \\
3       & 3 & 1 & 2 \\
\end{array}
& \ \ \ \ \ \ \ &
\begin{array}{c|ccc}
\circ & 1 & 2 & 3 \\ \hline\rule{0mm}{3.5mm}
1       & 1 & 2 & 3 \\
2       & 3 & 1 & 2 \\
3       & 2 & 3 & 1 \\
\end{array}
\end{array}
\]}}\noindent
The first has the following tracks: $\varphi_{1}=(1.23.)$,
$\varphi_{2}=(12.3.)$, $\varphi_{3}=(13.2.)$, the second:
$\psi_{1}=(1.2.3.)=\varepsilon$, $\psi_{2}=(123.)$,
$\psi_{3}=(132.)$. The first has no special tracks, the second has
one.
\end{example}

The above examples suggest that any unipotent quasigroup has a
special track. Indeed, if $x\cdot x=a$ for all $x\in Q$ and some
fixed $a\in Q$, then, as it is not difficult to see,
$\varphi_{a}=\varepsilon$ is its special track. Moreover,
$\varphi_{a}=\varepsilon$ if and only if $x\cdot x=a$ for all
$x\in Q$.
\begin{lemma}\label{lem}
If $\,\varphi_{k}$ is a special track of a quasigroup $Q(\cdot),$
then

$(a)$ \ $\alpha(k)=k$,

$(b)$ \ $\varphi_k\alpha=\alpha\varphi_k$,

$(c)$ \ $\varphi_k(k)=\alpha(\varphi_k(k))$

\noindent for any $\alpha\in {\rm Aut}\,Q(\cdot)$.
\end{lemma}
\begin{proof}
Indeed, $Z(\varphi_k)\ne Z(\varphi_i)=Z(\varphi_{\alpha(i)})$ for
every $i\ne k$ and $\alpha\in {\rm Aut}\,Q(\cdot)$ implies
$\alpha(i)\ne k$ for every $i\ne k$. Hence $\alpha(k)=k$. The
second statement is a consequence of \eqref{4}. $(c)$ follows from
$(a)$ and $(b)$.
\end{proof}

As a consequence of \eqref{4} and Lemma \ref{lem} $(a)$ we obtain
more general result.

\begin{proposition}\label{prop}
If $\,\varphi_{i}$, $\varphi_j$ are special track of a quasigroup
$Q(\cdot),$ then $$\varphi_i(j)=\alpha(\varphi_i(j))$$ for any
$\alpha\in {\rm Aut}\,Q(\cdot)$.\hfill $\Box{}$
\end{proposition}

\begin{example}\rm\label{ex2} The unipotent quasigroup from Example \ref{ex1} has no
special tracks. Its prolongation {\footnotesize
\[
\begin{array}{c|ccccc}
\cdot & 1 & 2 & 3 &4\\ \hline\rule{0mm}{3.5mm}
1       & 4 & 2 & 3&1 \\
2       & 3 & 1 & 4&2 \\
3       & 2 & 4 & 1&3 \\
4       & 1 & 3 & 2&4
\end{array}
\]}\noindent
obtained by the method proposed by Belousov (see \cite{bel} or
\cite{DD}) also has no special tracks.
\end{example}
\begin{example}\rm\label{ex3} The idempotent quasigroup of order $3$ has no special
track, but its prolongation obtained by Bruck's method (see
\cite{Br} or \cite{DD}) is an unipotent quasigroup with one
special track.
\end{example}
\begin{example}\rm\label{ex4} 
The cyclic group of order $4$ has no special tracks. Its
prolongation{\footnotesize
\[
\begin{array}{c|ccccc}
\cdot & 1 & 2 & 3 &4&5\\ \hline\rule{0mm}{3.4mm}
1       & 1 & 2 & 5&4&3 \\
2       & 2 & 3 & 4&5&1\\
3       & 3 & 4 & 1&2&5 \\
4       & 5 & 1 & 2&3&4\\
5       & 4 & 5 & 3&1&2
\end{array}
\]}\noindent
obtained according to the formula $(9)$ from \cite{DD} has three
special tracks:
\[
\varphi_{2}=(12.34.5.)\,, \ \ \ \varphi_{4}=(145.23.)\,, \ \ \
\varphi_{5}=(13524.)\,.
\]
\end{example}

\section{Rigid quasigroups}
\setcounter{theorem}{0}

Autotopies of a quasigroup form a group. Isotopic quasigroups have
isomorphic groups of autotopies (see for example \cite{bel} or
\cite{1}) but groups of automorphisms of such quasigroups may not
be isomorphic. Below we give examples of such quasigroups.

\begin{example}\rm\label{fig1} Let $Q(\cdot)$ be a quasigroup defined by the
following table:{\footnotesize
\[
\begin{array}{c|cccc}
\cdot & 1 & 2 & 3 & 4 \\ \hline\rule{0mm}{3mm}
1       & 1 & 2 & 3 & 4 \\
2       & 3 & 1 & 4 & 2 \\
3       & 4 & 3 & 2 & 1 \\
4       & 2 & 4 & 1 & 3 \\
\end{array}
\]}\noindent
It is not difficult to see that this quasigroup is isotopic to a
cyclic group of order $4$ and has the following four tracks:
\[
\varphi_{1}=(1.2.34.)\,, \ \ \ \varphi_{2}=(124.3.)\,, \ \ \
\varphi_{3}=(132.4.)\,, \ \ \ \varphi_{4}=(1423.)\,.
\]
Tracks $\varphi_{1}$ and $\varphi_{4}$ are special. So, according
to Lemma \ref{lem}, for any $\alpha\in{\rm Aut}\,Q(\cdot)$ we
have
$$
\alpha(1)=1\,, \hspace{5mm} \alpha(4)=4\,,
$$
which by Proposition \ref{prop} implies $\alpha(3)=3$. Hence
$\alpha(2)=2$, i.e., $\alpha=\varepsilon$. This means that this
quasigroup has only one (trivial) automorphisms while a cyclic
group of order $4$ has two automorphisms.
\end{example}

\begin{definition}\rm A quasigroup having only one automorphism is called
{\it rigid}.
\end{definition}

The above examples prove that a quasigroup isotopic to a rigid
quasigroup may not be rigid. Quasigroups of order two are rigid.

\begin{proposition}
No rigid quasigroups of order three.
\end{proposition}
\begin{proof} Indeed, if a quasigroup of order $3$ has an idempotent $e$ then
$\alpha=(e.xy.)$ is its non-trivial automorphism. If it has no
idempotents then it is commutative and has an automorphism
$\alpha=(123.)$.
\end{proof}

Each finite quasigroup containing at least $5$ elements is
isotopic to some rigid quasigroup \cite{Izb}. The same is true for
quasigroups defined on countable sets. So, {\it for every $k>3$
there exists at least one rigid quasigroup of order $k$}.

There are no rigid medial quasigroups of finite order $k\geqslant 2$
\cite{scerb}, but on the additive group of integers we can define
infinitely many rigid medial quasigroups \cite{soh}. A simple
example of such quasigroup is the quasigroup $(\mathbb{Z},\circ)$
with the operation $x\circ y=-x-y+1$. Finite rigid T-quasigroups
are characterized in \cite{scerb}.

Note, by the way, that prolongation does not save this property.
Nevertheless in some cases a prolongation of a rigid quasigroup
is also a rigid quasigroup. Moreover, a prolongation of a
non-rigid quasigroup may be a rigid quasigroup.

\begin{example}\label{ex5}\rm
The cyclic group of order $4$ is not a rigid quasigroup. Its
prolongation from Example \ref{ex4} is rigid. Indeed, it has three
special tracks $\varphi_2$, $\varphi_4$ and $\varphi_5$. Thus,
according to Lemma \ref{lem}, for any its automorphism $\alpha$ we
have $\alpha(2)=2$, $\alpha(4)=4$, $\alpha(5)=5$. Since
$\varphi_2(2)=1$, $\varphi_2(4)=3$, Proposition \ref{prop} implies
$\alpha(1)=1$ and $\alpha(3)=3$. Hence $\alpha=(1.2.3.4.5.)$,
which proves that this quasigroup is rigid.
\end{example}

\begin{example}\label{ex6}\rm
The loop $Q(\cdot)$ with the multiplication table {\footnotesize
\[
\begin{array}{c|cccccc}
\cdot & 1 & 2 & 3 & 4 &5&6\\ \hline\rule{0mm}{3mm}
1       &1&2&3&4&5&6 \\
2       &2&1&4&5&6&3 \\
3       &3&6&2&1&4&5 \\
4       &4&5&6&2&3&1 \\
5       &5&3&1&6&2&4 \\
6       &6&4&5&3&1&2
\end{array}
\]}
has the following tracks:
\[\begin{array}{llll}\varphi_{1}=(1.2.3465.), & \varphi_{2}=(12.3.4.5.6.),
&\varphi_{3}=(13.2645.), \\[3pt]
\varphi_{4}=(14.2356.), &\varphi_{5}=(15.24.36.),
&\varphi_{6}=(16.2543.)\,.
\end{array}
\]
Since
\[\begin{array}{llll}
Z(\varphi_{1})=[1,1,4],&Z(\varphi_{2})=[1,1,1,1,2],&Z(\varphi_{3})=[2,4],\\[3pt]
Z(\varphi_{4})=[2,4],&Z(\varphi_{5})=[2,2,2],&Z(\varphi_{6})=[2,4],
\end{array}
\]
tracks $\varphi_{1}$, $\varphi_{2}$, $\varphi_{5}$ are special.
So, according to Lemma \ref{lem}, for any automorphism  $\alpha$
of this quasigroup should be
\[
\alpha(1)=1, \hspace{5mm} \alpha(2)=2, \hspace{5mm} \alpha(5)=5.
\]
By Proposition \ref{prop}, we also have
$\alpha(3)=\alpha(\varphi_1(5))=\varphi_1(5)=3$ and
$\alpha(4)=\alpha(\varphi_1(3))=\varphi_1(3)=4$. Thus
$\alpha=(1.2.3.4.5.6.)=\varepsilon$, which means that this loop is
a rigid quasigroup.
\end{example}

In a similar way we can verify that the following four loops are
rigid: {\footnotesize
\[\begin{array}{lcr}
\begin{array}{c|cccccc}
\cdot & 1 & 2 & 3 & 4 &5&6\\ \hline\rule{0mm}{3mm}
1       &1&2&3&4&5&6 \\
2       &2&1&4&5&6&3 \\
3       &3&5&1&6&2&4 \\
4       &4&6&5&1&3&2 \\
5       &5&3&6&2&4&1 \\
6       &6&4&2&3&1&5
\end{array}
&&
\begin{array}{c|cccccc}
\circ & 1 & 2 & 3 & 4 &5&6\\ \hline\rule{0mm}{3mm}
1       &1&2&3&4&5&6 \\
2       &2&1&4&3&6&5 \\
3       &3&5&1&6&2&4 \\
4       &4&6&2&5&1&3 \\
5       &5&3&6&2&4&1 \\
6       &6&4&5&1&3&2
\end{array}
\\[15mm]
\begin{array}{c|cccccc}
*& 1 & 2 & 3 & 4 &5&6\\ \hline\rule{0mm}{3mm}
1       &1&2&3&4&5&6 \\
2       &2&3&6&1&4&5 \\
3       &3&4&5&2&6&1 \\
4       &4&5&2&6&1&3 \\
5       &5&6&1&3&2&4 \\
6       &6&1&4&5&3&2
\end{array}
&&
\begin{array}{c|cccccc}
\star & 1 & 2 & 3 & 4 &5&6\\ \hline\rule{0mm}{3mm}
1       &1&2&3&4&5&6 \\
2       &2&3&5&1&6&4 \\
3       &3&1&2&6&4&5 \\
4       &4&5&6&2&1&3 \\
5       &5&6&4&3&2&1 \\
6       &6&4&1&5&3&2
\end{array}\end{array}
\]}

We say that two quasigroups $Q(\cdot)$ and $Q(*)$ are {\it dual}
if
\[ x*y=y\cdot x\]
holds for all $x,y\in Q$. Dual quasigroups have the same
automorphisms. This means that {\it a quasigroup $Q(\cdot)$ is
rigid if and only if its dual quasigroup $Q(*)$ is rigid}.

\section{ Super rigid quasigroups}
\setcounter{theorem}{0}

The next interesting class of quasigroups is a class of
quasigroups having only one (trivial) autotopism. Quasigroups with
this property are called {\it super rigid}. Clearly, a super rigid
quasigroup has only one automorphism. Hence a super rigid
quasigroup is rigid. So, there are no super rigid quasigroups of order $2$
and $3$.

\smallskip
We remind some definitions and basic facts from \cite{3} and
\cite{2}.

\begin{definition}\rm By a {\it spin} of quasigroup $Q(\cdot)$ we mean the
permutation
$$
\varphi_{ij}=\varphi_{i}\varphi_{j}^{-1},
$$
where $\varphi_{i}$, $\varphi_{j}$ are tracks of $Q(\cdot)$. The
spin $\varphi_{ii}$ is called {\it trivial}.
\end{definition}

The set $\Phi_{Q}$ of all non-trivial spins of a quasigroup
$Q(\cdot)$ is called a {\it halo}. It can be divided into $n$
disjoint parts $\Phi_1,\Phi_2,\ldots,\Phi_n$, where
\[
\Phi_{i}=\{\varphi_{i1},\varphi_{i2},\ldots,\varphi_{i(i-1)},\varphi_{i(i+1)},
\ldots,\varphi_{in}\}.
\]

Let $\Phi=\{\sigma_{1},\sigma_{2},\ldots,\sigma_{k}\}\subseteq
S_Q$ be a collection of permutations of the set $Q$. According to
\cite{2}, the set
\[
Sp\,(\Phi)=[Z(\sigma_{1}),Z(\sigma_{2}),\ldots,Z(\sigma_{k})]
\]
is called the {\it spectrum} of $\Phi$. The spectrum of all spins
of $Q(\cdot)$ is called the {\it spin-spectrum}.

\begin{example}\label{ex7}\rm The quasigroup considered in the
Example \ref{fig1} has the following proper spins:
\[\begin{array}{ll}
\varphi_{12}=(1342.), \ \ \varphi_{13}=(1243.), \ \
\varphi_{14}=(14.23.), \\[3pt]
\varphi_{21}=(1243.), \ \ \varphi_{23}=(14.23.), \ \
\varphi_{24}=(1243.), \\[3pt]
\varphi_{31}=(1342.), \ \ \varphi_{32}=(14.23.), \ \
\varphi_{34}=(1243.),\\[3pt]
\varphi_{41}=(14.23.), \ \ \varphi_{42}=(1243.), \ \
\varphi_{43}=(1342.). \end{array}
\]
Thus\ $Sp\,(\Phi_1)=[[4],[4],[2,2]]$,
$Sp\,(\Phi_2)=[[4],[2,2],[4]]=Sp\,(\Phi_3)$,
$Sp\,(\Phi_4)=[[2,2],[4],[4]]$. In the abbreviated form it will be
written as $Sp\,(\Phi_i)=2A+B$, where $A=[4]$, $B=[2,2]$.
\end{example}

Finite isotopic quasigroups have the same spin-spectrum (\cite{2},
Theorem 2.5). Moreover, spins of isotopic quasigroups are pairwise
conjugated. Namely, if quasigroups $Q(\cdot)$ and $Q(\circ)$ are
isotopic and
$$
\gamma(x\circ y)=\alpha(x)\cdot\beta(y),
$$
then spins $\varphi_{ij}$ of $Q(\cdot)$ and $\psi_{ij}$ of
$Q(\circ)$ are connected by the equality
\[
\varphi_{\gamma(i)\gamma(j)}=\beta\psi_{ij}\beta^{-1}.
\]
Hence, identifying $Q(\cdot)$ and $Q(\circ)$, we obtain
\begin{equation}\label{spin}
\varphi_{\gamma(i)\gamma(j)}=\beta\varphi_{ij}\beta^{-1}.
\end{equation}
This means that for any fixed $i\in Q$ and an arbitrary
permutation $\gamma$ of $Q$, we have
\[
Sp\,(\Phi_{i})=Sp\,(\Phi_{\gamma(i)}).
\]

If\ $Sp\,(\Phi_i)\neq Sp\,(\Phi_k)$ for all $k\in Q$, $k\neq i$,
then we say that the part $\Phi_i$ is {\it special}.

It is not difficult to see that the following lemma is valid.

\begin{lemma}\label{lem2} If $\,\Phi_{i}$ is a special part of\
$\Phi_Q$, then $\gamma(i)=i$ for any autotopism $(\alpha,\beta,
\gamma)$ of $\,Q(\cdot)$. \hfill $\Box{}$
\end{lemma}

\begin{proposition}\label{P-dual}
Dual quasigroups have the same spin-spectrum and their special parts
have the same numbers.
\end{proposition}
\begin{proof}
Let $Q(\cdot)$ and $Q(\circ)$ be dual quasigroups. If $\varphi_i$
is a track of $Q(\cdot)$, then
$$
\varphi_i(x)\circ x=x\cdot\varphi_i(x)=i
$$
for every $x\in Q$. From this, replacing $x$ by
$\varphi_i^{-1}(x)$, we obtain $x\circ\varphi_i^{-1}(x)=i$, which
means that $\psi_i=\varphi_i^{-1}$ is a track of $Q(\cdot)$. So,
spins $\psi_{ij}$ of $Q(\circ)$ have the form
$$
\psi_{ij}=\psi_i\psi_j^{-1}=\varphi_i^{-1}\varphi_j=(\varphi_j^{-1}\varphi_i)^{-1}.
$$

Since for any conjugate permutations $\sigma_1$, $\sigma_2$ of the
same set $Q$ we have $Z(\sigma_1)=Z(\sigma_2)$ (cf. \cite{4}), for
any permutations $\alpha$, $\beta$, from
$\alpha\beta=\beta^{-1}(\beta\alpha)\beta$ it follows
$Z(\alpha\beta)=Z(\beta\alpha).$ Thus
$$
Z(\psi_{ij})=Z((\varphi_j^{-1}\varphi_i)^{-1})=Z(\varphi_j^{-1}\varphi_i)=Z(\varphi_i\varphi_j^{-1})
=Z(\varphi_{ij}),
$$
for $i,j=1,2,\ldots,n$. Consequently $Sp\,(\Psi_i)=Sp\,(\Phi_i)$
for all $i=1,2,\ldots,n$.
\end{proof}

\begin{proposition}\label{dual}
A quasigroup $Q(\cdot)$ is super rigid if and only if its dual
quasigroup $Q(\circ)$ is super rigid.
\end{proposition}
\begin{proof}
Let $Q(\cdot)$ be a super rigid quasigroup. If
$(\alpha,\beta,\gamma)$ is an autotopism of a dual quasigroup
$Q(\circ)$, then $(\beta,\alpha,\gamma)$ is an autotopism of
$Q(\cdot)$. Hence $\alpha=\beta=\gamma=\varepsilon$.
\end{proof}

Now we give examples of super rigid quasigroups.

\begin{example}\label{ex8}\rm
Consider the following quasigroup: {\footnotesize \[
\begin{array}{c|ccccccc}
\cdot & 1 & 2 & 3 & 4 & 5 & 6 & 7\\ \hline\rule{0mm}{3mm}
1       & 1 & 2 & 3 & 4 & 5 & 6 & 7\\
2       & 2 & 1 & 7 & 6 & 4 & 5 & 3\\
3       & 3 & 6 & 1 & 2 & 7 & 4 & 5\\
4       & 4 & 5 & 2 & 1 & 3 & 7 & 6\\
5       & 5 & 7 & 4 & 3 & 6 & 2 & 1\\
6       & 6 & 3 & 5 & 7 & 2 & 1 & 4\\
7       & 7 & 4 & 6 & 5 & 1 & 3 & 2
\end{array}
\]}
This quasigroup has seven tracks:
\[
\begin{array}{lllll}
\varphi_{1}=(1.2.3.4.57.6.), & \varphi_{2}=(12.34.56.7.), &
\varphi_{3}=(13.276.45.), \\[3pt]
\varphi_{4}=(14.25367.), & \varphi_{5}=(15.26374.),&
\varphi_{6}=(16.2473.5.), \\[3pt] \varphi_{7}=(17.235.46.).
\end{array}
\]

After the calculation of all spins we can see that each spin can
be decomposed into cycles in one of the following ways:
$$
A=[7],\hspace{1cm}B=[3,4],\hspace{1cm}C=[2,2,3],\hspace{1cm}D=[2,5].
$$
Moreover,
$$\arraycolsep.5mm
\begin{array}{lrcrcrcr}
Sp\,(\Phi_{1})= & A  &  &  & +&2C & +&3D, \\[2pt]
Sp\,(\Phi_{2})= & A  & +&B & +&2C & +&2D, \\[2pt]
Sp\,(\Phi_{3})= &    & &2B & +&2C & +&2D, \\[2pt]
Sp\,(\Phi_{4})= & 2A & &   & +&C  & +&3D, \\[2pt]
Sp\,(\Phi_{5})= & 2A &+&B & +&C  & +&2D, \\[2pt]
Sp\,(\Phi_{6})= &    & &   & &2C & +&4D, \\[2pt]
Sp\,(\Phi_{7})= &    & &   & &2C & +&4D.
\end{array}
$$
Since parts $\Phi_{1},\Phi_{2},\Phi_{3},\Phi_{4},\Phi_{5}$ are
special, from Lemma \ref{lem2} it follows that for any autotopism
$(\alpha,\beta,\gamma)$ of this quasigroup we have
$\gamma=(1.2.3.4.5.67)$ or $\gamma=\varepsilon$.

Below we prove that the first case is impossible. For this we
consider two spins
$$
\varphi_{15}=(1736245.) \ \ \ {\rm and } \ \ \
\varphi_{52}=(16.153.47.).
$$
According to \eqref{spin}, we have
$$
\varphi_{15}=\beta\varphi_{15}\beta^{-1} \ \ \ {\rm and } \ \ \
\varphi_{52}=\beta\varphi_{52}\beta^{-1}.
$$
In view of Theorem 5.1.3 from \cite{4} any $\beta$ satisfying the
first equality has the form
\[
\beta=\varphi_{15}^{i},\hspace{5mm}i=1,2,3,\ldots,7.
\]
The second equality is satisfied by $\beta=\varphi_{52}^{j}$. So,
$\varphi_{15}^{i}=\varphi_{52}^{j}$ for some $i,j$. Since
$\varphi_{52}^{j}(1)=6$ or $\varphi_{52}^{j}(1)=1$, we have
$\varphi_{15}^{i}(1)=6$ or $\varphi_{15}^{i}(1)=1$. The first case
holds for $i=3$, the second -- for $i=7$. The case $i=3$ is
impossible because $\varphi_{15}^{3}(6)=5\ne\varphi_{52}^{j}(6)$.
So, $i=7$ and $\beta=\varphi_{15}^{7}=\varepsilon$. Thus
\[
\gamma(x\cdot y)=\alpha(x)\cdot y,
\]
which implies $\gamma(x)=\gamma(x\cdot 1)=\alpha(x)$ for every
$x\in Q$. Consequently,
\[
\gamma(6)=\gamma(3\cdot 2)=\alpha(3)\cdot 2=\gamma(3)\cdot
2=3\cdot 2=6.
\]
Hence $\gamma=\alpha=\varepsilon$. This proves that this
quasigroup is super rigid.
\end{example}

It is the smallest rigid quasigroup. To prove this fact first we
select all rigid quasigroups of order $k<7$, next we prove that
these quasigroups are not super rigid.

\begin{example}\label{ex9}\rm Consider the quasigroup:
 {\footnotesize \[
\begin{array}{c|cccccccccccc}
\cdot &1&2&3&4&5&6&7&8&9\\ \hline\rule{0mm}{3mm}
1     &1&2&3&4&5&6&7&8&9\\
2     &2&3&1&8&6&7&5&9&4\\
3     &3&1&2&9&7&5&6&4&8\\
4     &4&5&6&7&9&8&1&3&2\\
5     &5&6&4&2&1&9&8&7&3\\
6     &6&4&5&3&8&1&9&2&7\\
7     &7&8&9&5&3&2&4&6&1\\
8     &8&9&7&1&4&3&2&5&6\\
9     &9&7&8&6&2&4&3&1&5
\end{array}
\]}

Using the same method as in Example \ref{ex8} we can see that
$Sp\,(\Phi_3)= Sp\,(\Phi_4)$ and $Sp\,(\Phi_i)\ne Sp\,(\Phi_j)$
for all  $i\ne j$, \ $i\ne 3,4$. This means that
$\Phi_1,\Phi_2,\Phi_5,\Phi_6,\Phi_7,\Phi_8$ and $\Phi_9$ are
special. Thus, by Lemma \ref{lem2}, for any autotopism
$(\alpha,\beta,\gamma)$ of this quasigroup should be
$\gamma=(1.2.34.5.6.7.8.9.)$ or $\gamma=\varepsilon$.

We prove that $\gamma=\varepsilon$. For this consider two spins
$$
\varphi_{68}=(197286345.) \ \ \ {\rm and } \ \ \
\varphi_{13}=(123.46.59.78.).
$$
Then, similarly as in the previous example,
$\varphi_{68}=\beta\varphi_{68}\beta^{-1}$ and
$\varphi_{13}=\beta\varphi_{13}\beta^{-1}$ imply
\[
\beta=\varphi_{68}^{i}=\varphi_{13}^{j}
\]
for some $i=1,2,3,\ldots,9$ and $j=1,2,\ldots,6$. Since
$\varphi_{13}^{j}(4)=4$ or $\varphi_{13}^{j}(4)=6$, also
$\varphi_{68}^{i}(4)=4$ or $\varphi_{68}^{i}(4)=6$. Thus $i=9$ or
$i=7$. For $i=7$ we have $\varphi_{68}^7(3)=8$. But
$\varphi_{13}^j(3)\ne 8$ for every $j$. So, this case is
impossible. Therefore $i=9$. Consequently
$\beta=\varphi_{68}^9=\varepsilon$, i.e., $\gamma(x\cdot
y)=\alpha(x)\cdot y$, which implies $\gamma(x)=\gamma(x\cdot
1)=\alpha(x)$ for every $x\in Q$. Now, using the above we obtain
$$
\gamma(3)=\gamma(2\cdot 2)=\alpha(2)\cdot 2=\gamma(2)\cdot
2=2\cdot 2=3.
$$ Hence $\gamma=\alpha=\varepsilon$. This means
that this quasigroup has no nontrivial autotopisms. So, it is
super rigid.
\end{example}

\noindent \footnotesize{\rightline{Received \ March 28, 2009}

\medskip\noindent
A.I.Deriyenko and I.I.Deriyenko:\\ Department of Higher
Mathematics and Informatics, Kremenchuk State Polytechnic
University,
20 Pervomayskaya str, 39600 Kremenchuk, Ukraine\\
E-mails: andrey.deriyenko@gmail.com, \ \ ivan.deriyenko@gmail.com\\[5mm]
W.A.Dudek:\\ Institute of Mathematics and Computer Science,
Wroclaw University of Technology, Wyb. Wyspia\'nskiego 27, 50-370
Wroclaw, Poland// E-mail: dudek@im.pwr.wropc.pl }

\end{document}